\font\eightrm=cmr8         
\font\eightit=cmti8        
\font\eightmi=cmmi8        


\def\frac#1#2{{#1\over #2}}   

\def\N{{\bf N}}               

\def\qed{\allowbreak\qquad\null
         \nobreak\hfill\square}
\def\R{{\bf R}}               
\def\square{\vrule height 5pt depth 0pt
             width 5pt}       

\def\Tn{{\cal T}}             
\def\.{\cdot}                 
\def\:{\colon}                
\def\<#1>{\langle#1\rangle}   

\def\then{\Rightarrow}

\def\SDS{{\eightrm TDS}}

\outer\def\section#1. #2\par{
      \bigskip\bigskip       
      \message{#1. #2}
      \leftline{\bf#1. #2}
      \nobreak\smallskip     
      \noindent}             

\def\declare#1. #2\par{
      \medskip\noindent     
      {\bf#1.}\rm           
      \enspace\ignorespaces 
      #2\par\smallskip}     

\def\cite#1{{\rm[#1]}}        
\def\refno#1. #2\par{\smallskip\item{[#1]} #2\par}


\def\Bastiani{1}
\def\Boman{2}
\def\Choquet{3}
\def\DonatoGeometrie{4}
\def\DonatoVirasoro{5}

\def\Fro{7}
\def\Frolicherb{8}
\def\Bfroli{9}
\def\Hamilton{10}
\def\Jarchow{11}
\def\Ke{12}
\def\Keller{13}
\def\Ckeller{14}
\def\Kri{15}
\def\Krieglb{16}
\def\Lang{17}
\def\Lawvere{18}
\def\Marinescu{19}
\def\Milnor{20}
\def\Seip{21}
\def\SouriauAlgoritme{22}
\def\Sou{23}
\def\Torre{24}
\def\TorreTangen{25}
\def\TorreBan{26}
\def\Wong{27}
\def\Yamamuro{28}

\newtoks\rightheadtext      
\newtoks\leftheadtext       

\headline={\ifnum\pageno>1  
            \ifodd\pageno   
             \hfil \eightrm \the\rightheadtext \hfil
             \llap{\tenrm\folio}%
            \else \rlap{\tenrm\folio}
             \hfil \eightrm \the\leftheadtext \hfil \fi
           \else\hfil \fi}
\footline={\hfil}           

\rightheadtext{Differentiability, Convenient Spaces}

\leftheadtext{CARLOS A. TORRE}



\centerline{\bf Differentiability, Convenient Spaces and Smooth Diffeologies }

\bigskip

\centerline{\sl Carlos A. Torre}
\bigskip
\centerline{\sl Escuela de Matem\'atica,
                 Universidad de Costa Rica,}
\centerline{\sl San Jos\'e, Costa Rica}
\bigskip

{\narrower
\eightrm
\baselineskip=9.5pt
\textfont1=\eightmi
\noindent
{\eightit ABSTRACT}.

We review the basic definitions and properties concerning smooth
structures, convenient spaces, diffeological spaces and tangent
structures. The relation between the first two is described. A tangent
structure is constructed for each pre-convenient space. This one is
proved to be convenient if and only if the space and the tangent fibres
coincide.

\par
\smallskip}

\bigskip

\section 1. Introduction

Infinite dimensional spaces appear in the description of many physical
systems as configuration spaces, symmetry or gauge groups. However,
most of these function spaces are not Banach and 
unfortunately, classical differentiation uses Banach spaces.
The difficulty arises because 
the classical condition of Fr'chet differentiability
requires an appropiate structure on the space $L(E,F)$ of linear
morphisms, which does not exist if the spaces are not Banach
(a straightforward generalization is not possible).
However, geometrical difficulties arise in the category of
manifolds which are not Banach, including the study of submanifolds,
orbit spaces, exponential map, etc.

For this reason, many attempts, successful in different ways, 
have been made to substitute the Banach space theory of
differentiation and manifolds, at different levels of generalization.
In the area of differentiation it includes topology, convergence
structures, bornologies, etc. Other approaches, more related with 
spaces of differentiable functions, 
include smooth structures (and convenient
spaces) and (tangent) diffeological spaces. 

Convenient spaces \cite{\Fro, \Krieglb, \Seip} 
are in many ways one of the most successful and natural
generalizations of Banach spaces. The basic idea is based on the
following two facts:
(a) a curve in a Banach space $c\colon \R \to E$ is smooth iff
$l\circ c$ is smooth iff $\forall l\in E^{\prime}$, (b) a map between
two Banach spaces $g\colon E \to F$ is smooth iff $g\circ c$ is smooth for
every smooth curve $c\colon \R \to E$. This indicates that only
$E^{\prime}, F^{\prime}$ are required, without the use of norms and
topologies. Diffeological spaces \cite{\SouriauAlgoritme, 
\DonatoGeometrie, \Sou} instead, were defined in a more geometric
setting, with the purpose of obtaining geometric quantization in
infinite dimensional cases. In these spaces, the construction of
geometric objets is based on a set of plaques, generalizing the basic
constructions in finite dimensional manifolds (with the use of charts).

Our purpose here is to study some general and basic relations between these
two approaches.
In section 2 we recall the basic properties of the differentiable
functions $f\colon E \to F$ in locally convex spaces 
with respect to several locally
convex topologies and convergence structures in the space 
${\cal L}(E, F)$ of continuous
linear maps between these spaces. Section 3 describes the basic
definitions and properties of smooth structures. The definition of
a diffeological space is recalled in section 4, and section 5 describes
a relation between the these two categories. Section 6 defines a
tangent structure on a diffeological space. In section 7 we review
the definition and basic properties of convenient spaces. Section 8
describes a relation between convenient spaces and tangent structures.
The next section describes a relation between several definitions of
a tangent vector in these categories and the basic relations between
manifolds and tangent diffeological spaces.

\section 2. Differentiation in locally convex spaces.

The definition of the derivative of a function 
$f\colon U\subset \R \to \R$, $U$ open, at a point $t_0$:
$\lim_{t\to t_0}{f(t_0+t)-f(t_0)\over t-t_0}$ is generalized to a function
$f\colon U\subset \R \to F$ if $F$ is a topological vector space, by
changing the topology of the range in the limit and if the domain is
a topological vector space $E$, then the directional 
derivative is defined along each vector $v\in E$: 
$df(x, v) = (f\circ c)^{\prime}(0)$
where $c(t) = x+tv$. It is better to require $E, F$ to be Hausdorff 
in order to have uniqueness of limits and 
also they are required to be locally convex vector spaces (l.c.s.), 
otherwise
there may be functions whose derivatives are cero everywhere but that are
not constant. Inductively, if $k\in \N$ and $v_i\in E, i=1,\cdots k$ then
$$d^kf(x,v_1, \cdots ,v_k)= d(d^{k-1}f(x,v_1, \cdots ,v_{k-1}))(v_k)$$

The existence of all directional derivatives is a weak condition
since it does not even imply continuity
of the function (consider the function $f(0,0)=0$ and
$f(x,y)= {xy^2\over x^2+y^4}$ otherwise).

The map $v\to df(x,v)$ is not linear in general if $df(x,v)$
exists for all $v\in E$, not even in the finite dimensional case
(consider the function $f(0,0)= 0$ and $f(x,y)= {xy^2\over x^2+y^2}$
otherwise).

Besides the existence of $d^kf(x, v_1, \cdots v_k)$ for all $x\in U$
and $v_i\in E$, we usually require a continuity condition: $f$ is
called $C^k$-{\it Gateaux differentiable} if in addition the map
$d^kf\colon U\times E^n \to F$
is continuous.
In this case $d^kf$ is multilinear and symmetric and the
space of $C^k$-maps is closed under composition
\cite\Milnor and  chain rule holds. Inverse function
theorem holds in the category of tamed Fr'chet spaces \cite\Hamilton.

In the case of Banach
spaces a stronger condition of differentiability may be used:  denote,
in the general case, by $L^n(E,F)$ (resp. ${\cal L}(E,F)$) the vector
space of n-linear mappings (resp. continuous) from $E$ to $F$. Then
the map $f$ is
called {\it Fr'chet differentiable} at $x_0\in U$ if
there exists a linear continuous map (equal to $df(x_0, \cdot)$) which is
tangent to $f(x_0+v)-f(x_0)$ at $v=0$ (this is the
remainder condition) and it is called Fr'chet $C^1$ 
on $U$ if the map
$x_0\to df(x_0, \cdot )\in L(E, F)$ is continuous (when $L(E, F)$ is
given the canonical Banach space structure), denoted $Df(x_0)$.

The second derivative is a map
$x_0 \to L(E, L(E, F))$ and inductively, the map $D^kf\colon U \to L(E^k, F)$
is required to be continuous. 

An important fact is that the remainder condition is a consecuense of
continuity of $Df$ (this simplifies the study of generalizations to 
the study of different conditions of continuity of $Df$ only).

Composition of $C^k$ maps is again
$C^k$, chain rule, implicit function theorem and inverse function
theorem holds, similar to the finite dimensional case \cite\Lang.
Also good existence and uniqueness theorem of ordinary differential
equations holds \cite\Choquet.

Fr'chet differentiability at a point implies the existence of all
directional derivatives and continuity of $f$ but does not imply continuity
of directional derivatives (even if $f$ is Fr'chet differentiable in an
open set (take $f(0,0)=0$ and $f(x,y)= {x^2y^2\over x^2+y^2}$ otherwise).
A weak converse is well known in the finite dimensional case (if all partial
derivatives of $f$ exist and are continuous in an open set $U$ then
$f$ is Fr'chet differentiable at each point of $U$). However, existence
of all directional derivatives and continuity of them along lines
through a point does not imply Fr'chet differentiability at the point.

In the finite dimensional case Gateaux $C^k$
and Frechet $C^k$ are equivalent. In the 
more general case of Banach spaces both are equivalent if the condition
`$df(x_0, v)$ continuous' is replaced
by `$df(x_0, v)$ $Lip^k$'(this means for normed spaces that in the 
definition of Gateaux differentiability the condition of continuous
is replaced by Lipschitz locally at each point). 

The notion of Fr'chet differentiability may not be generalized,
not even to Fr'chet spaces, because the dual of a Fr'chet space which
is not Banach is never a Fr'chet space. 

For functions between Banach spaces, the continuity of
$df\colon E\times E \to F$ does not imply continuity of
$Df\colon E \to L(E,F)$ (however $Df$ is locally Lipschitz iff
$df$ is locally Lipschitz).
For this reason an alternative continuity condition 
for any $E, F$ l.c.s. is constructed in the following way:
$f$ is called
{\it weakly p-times differentiable on U} if 
$$D^{k+1}f(x)v := \lim_{t\to 0}t^{-1}[D^kf(x+tv)-D^kf(x)]$$ exists
$\forall (x,v)\in U\times E$ and $k= 0, \cdots ,p-1$, with respect to 
the topology of simple convergence in $L^k(E,F)$ and
$D^kf(x)v$ is linear with respect to $v\in E$, $\forall x\in U$.

A stronger condition 
is defined by providing a
l.c.s. topology on ${\cal L}(E,F)$: a collection $\sigma$
of bounded subsets of $E$ defines a family of seminorms 
on ${\cal L}(E, F)$ (the
supremum, for each seminorm of $F$, over subsets $S\in \sigma$). The 
space with the topology $\sigma$ is denoted 
${\cal L}_{\sigma}^n(E,F)$. The most
important examples are given by finite, compact, precompact, and bounded
sets, whose topology is denoted by $\wedge_s, \wedge_k, 
\wedge_{pk}, \wedge_b$. Also we define the spaces 
${\cal H}_{\sigma}^0(E,F):=F$ and 
$${\cal H}^n_{\sigma}(E, F):= 
{\cal L}_{\sigma}(E, {\cal H}^{n-1}_{\sigma}(E,F))$$
If $E$ is metrizable and barrelled, or if $E$ is metrizable and
$\wedge_k\subset \sigma$ then 
${\cal L}_{\sigma}^n(E, F)= {\cal H}_{\sigma}^n(E, F)$
We may also define a convergence structure $\wedge$ on ${\cal L}^n(E,F)$,
even if $E,F$ are convergence vector spaces.
Some examples are the continuous convergence $\wedge_c$ of
Bastiani \cite\Bastiani (this is the coarsest convergence
vector space structure (c.v.s.) for which
the evaluation is continuous), the quasi-bounded
convergence $\wedge_{qb}$ \cite\Bfroli , 
the bounded convergence $\wedge_b$ \cite\Yamamuro (which
coincides with the topology of bounded topology if $E, F$ are l.c.s.),
the Marinescu convergence $\wedge_{\triangle}$ \cite\Marinescu.
And the space with this structure is denoted ${\cal L}^n_{\wedge}(E,F)$.
We have $\wedge_c \leq \wedge_{qb} \leq \wedge_{\pi} \leq \wedge_{\triangle}$.
For these $c.v.s.$ and also if $E$ is Banach, we have 
${\cal L}^m(E, {\cal L}_{\wedge}^n(E,F)) \simeq {\cal L}^{m+n}(E,F)$.

Now, having a topology or convergence structure $\wedge$, a map 
$f\colon U\subset E \to F$ is called 
{\it differentiable of class $C^p_{\wedge}$}
if $f$ is weakly p-times differentiable and 
$D^kf(x)\in {\cal H}_{\wedge}^k(E,F)$, $\forall x\in U$ and $0\leq k\leq p$,
and $D^kf\colon U \to {\cal H}^k_{\wedge}(E,F)$ is continuous. 
If $\wedge \geq \wedge_c$ we may take ${\cal L}_{\wedge}^k$ instead.

We have the inclusions \cite\Ke \ : 
$$C^p_{\triangle} \subset C^p_{\pi} \subset C^p_{qb} \subset \{ C^p_c, C^p_b\} 
\subset C^p_b\subset C^p_{pk} \subset C^p_k \subset C^p_s$$
If $E$ is normable then 
$C^p_{\triangle}= C^p_{\pi} = C^p_{qb} = C^p_b$ and
$C^p_{pk}= C^p_k$. Also $C^p_{pk}=C^p_k=C^p_s$
if $E$ is Fr'chet; however, $C^1_s \neq C^1_b$ even if $E, F$ are Banach
(all of them are equivalent however if $E$ is finite dimensional). Therefore
we have two notions of differentiability in Fr'chet spaces: Gateaux-Levi 
($C^p_s$)
and Fr'chet ($C^p_b$), which coincide if $E$ is finite dimensional or
if $E$ is Fr'chet-Schwartz and $F$ is normable \cite\Wong.

Differentiability in l.c.s. preserves some properties: if $f\in C^p_s$
then $D^kf(x)$ is symmetric for all $0\leq k\leq p$ and Taylor approximation
holds. Also $f\in C^p_{\sigma}$ implies that $D^kf(x)v$ exists uniformly
with respect to $v\in S\in \sigma$ and it is linear with respect to $f$.
In addition, the remainder satisfies:  if $f\in C^p_{\sigma}$ then
$$\lim_{t\to 0}t^{-1}[D^kf(x+tv)-D^kf(x)-D^{k+1}f(x)tv]=0$$
uniformly with respect to $v\in S\in \sigma , \forall k\leq p-1$.

There are many limitations however: the evaluation map 
$${\cal L}^n_{\sigma}(E, F)\times E^n \to F$$ is not continuous for any 
l.c.s. topology on ${\cal L}^n(E,F)$ unless $E$ is normable \cite\Keller
(it is
continuous however for convergence structures stronger or equal to 
$\wedge_c$). Composition 
$${\cal L}_{\sigma}(F,G)\times {\cal L}_{\sigma}(E,F) \to {\cal L}_{\sigma}
(E,G)$$ 
is not continuous for any l.c.s. topology on ${\cal L}(E,F)$,
except for special $E,F$, such as Banach spaces \cite\Ckeller \
(it is continuous for 
convergence structures $\wedge_{qb}, \wedge_{\triangle}$). Differentiable
functions with respect to convergence structures stronger or equal to 
$\wedge_c$ are closed under composition, but it does not hold in the
category of l.c.s. unless they are Banach. Some restricted cases hold:
if $E$ is metrizable, $f\colon U\to F$ is $C^1_s$ and 
$g\colon V\subset F\to G$
is $C^1_k$ then $g\circ f\in C^1_s$. In addition the continuity of
$D^kf\colon U \to {\cal L}_{\wedge}(E,F)$ does not imply the continuity
of $d^kf\colon U\times E^n \to F$ in general. It holds 
however if $E$ is metrizable and barrelled and $\wedge= \wedge_s$ or if
$E$ is metrizable and $\wedge = \wedge_k$ or for general l.c.s. and 
$\wedge$ finer than $\wedge_c$.
Also the continuity of
$D^pf\colon U\to {\cal L}^p_s(E,F)$ does not imply the continuity of
$D^kf$ if $k<p$ (it does if $E$ is finite dimensional).

A map $f\colon U\to F$ is called of class $C^{\infty}_{\wedge}$ if
it is of class $C^p_{\wedge}$ for all $p\in \N$. The inclusions,
described above,
between spaces $C^p_{\wedge}$ for different $\wedge$, also holds
for $p=\infty$. 
If $E$ is Banach, $F$ l.c.s. or $E$ Fr'chet and $F$ normable then all
spaces $C^{\infty}_{\wedge}$ coincide, escept for
$C_{\triangle}^{\infty}$ and therefore we speak of the
space $C^{\infty}$ instead, which is closed under composition.
  In addition the composition of 
two $C^{\infty}_s$ functions is a $C^{\infty}_s$ function in the
category of Fr'chet spaces. 
The spaces are closed under composition in other special cases: if 
$f\in C^{\infty}_{\wedge}(E,F)$ and $g\in C^{\infty}_{\wedge}(F,G)$
then $g\circ f\in C^{\infty}_{\wedge}(E,G)$ if
$\wedge = \wedge_c$.

\section 3. Smooth structures.

In this section we recall some basic definitions and results concerning
smooth structures: a theory that provides a generalization
to the construction of $Lip^k$ and $C^{\infty}$ curves in Banach spaces.

In \cite\Boman \ Boman proves that a function $\R^n \to \R$ is 
smooth if and only if it is smooth along smooth curves. A result that
is generalized  \cite\Fro \
for any classical smooth manifold or any Banach space $X$: the family of
smooth curves $c\colon \R \to X$ 
and the family of smooth functions $f\colon X \to \R$ determine each other
(a map belongs to one of the families if the composites with the members
of the other family are functions in $C^{\infty}(\R, \R)$). Another
generalization of Boman's result is the following: 
a map between two Banach spaces
is a smooth function if and only if its composite with smooth curves
of the source and the smooth real valued functions of the range are
functions in $C^{\infty}(\R, \R)$.

These two results, together with the need of a differentiation theory in
general spaces that include spaces of functions such as diffeomorphisms
(which in general are not Banach) is the reason for a 
generalization: instead of
$C^{\infty}(\R, \R)$, a general set $\cal M$ of functions from a 
set $S$ to another set $R$ is considered, a Banach space is substituted
by a general set $X$, together with a family ${\cal C}_X$ of 
{\it curves}
$c\colon S \to X$ and a family ${\cal F}_X$ of 
{\it functions} $f\colon X \to R$ such that
one of the families determine the other one: given $X$, $\cal M$ we
observe that

\item{1.}
A set $\cal C$ of curves in $X$ determines a set $\cal F$ of functions
on $X$: 
$\{ f\colon X \to \R / f\circ c\in {\cal M}, \forall
c\in {\cal C} \}$, denoted $\Phi {\cal C}$.

\item{2.}
A set $\cal F$ of functions on $X$ determines a set $\cal C$ of curves in
$X$:
$\{ c\colon I \to X / f\circ c\in {\cal M}, \forall
f\in {\cal F} \}$, denoted $\Gamma {\cal F}$.

An $\cal M$-{\it structure} on $X$ \cite\Fro , \cite\Kri \ is a pair 
$({\cal C}, {\cal F})$ of
curves $\cal C$ in $X$ and functions $\cal F$ on $X$ such that they
determine each other: $\Phi {\cal C} = {\cal F}$ and
$\Gamma {\cal F} = {\cal C}$. Also $(X, {\cal C}, {\cal F})$ is called
an $\cal M$-space. A map $f\colon X \to Y$ between two spaces $X, Y$
with $\cal M$-structures $({\cal C}_X, {\cal F}_X)$ and
$({\cal C}_Y, {\cal F}_Y)$ respectively,
is called an $\cal M$-{\it map} if $h_*({\cal C}_X) \subset {\cal C}_Y$
(this is equivalent to $h^*({\cal F}_Y)\subset {\cal F}_X$ and also 
is equivalent to
$g\circ h\circ c \in {\cal M}, \forall c\in {\cal C}_X, \forall g\in
{\cal F}_Y$). The set of ${\cal M}$-maps $f\colon X \to Y$ is
denoted ${\cal M}(X, Y)$.
The elements of $\cal C$ are called {\it structure curves} and the
elements of $\cal F$ are called {\it structure functions}.

We may have several ${\cal M}$-structures on the same set $X$. In this
case a structure $({\cal C}_1, {\cal F}_1)$ on $X$ is called smaller
than $({\cal C}_2, {\cal F}_2)$ if ${\cal C}_1 \subset {\cal C}_2$.

For any set ${\cal C}_0$ of curves on $X$ there exists a smallest 
$\cal M$-structure $({\cal C}, {\cal F})$ on $X$ such that
${\cal C}_0 \subset {\cal C}$ (given by ${\cal F} = \Phi {\cal C}_0$
and ${\cal C} = \Gamma {\cal F}$). With this $\cal M$-structure any
map $h\colon X \to Y$ is an $\cal M$-map if and only if
$h_*({\cal C}_0) \subset {\cal C}_Y$.

For any set ${\cal F}_0$ of functions on $X$ there is a largest 
${\cal M}$-structure
$({\cal C}, {\cal F})$ such that ${\cal F}_0 \subset {\cal F}$ (given by
${\cal C} = \Gamma {\cal F}_0$ and ${\cal F} = \Phi {\cal C}$).
With this $\cal M$-structure a map $h\colon Y \to X$ is an $\cal M$-map
if and only if $h^*({\cal F}_0) \subset {\cal F}_Y$.
In particular if $X$ is a vector space and ${\cal F}_0$ is a set of
linear functions, we say that $(X, {\cal C}, {\cal F})$ is {\it linearly
generated} and the ${\cal M}$-space 
may be denoted by $(X, {\cal F}_0)$ instead.

Some fundamental examples of sets $\cal M$ related with differentiability are:
\item{(a)}
$C^k(\R, \R)$, denoted ${\cal M}^k$, for any $k\in \N \cup \{ \infty \}$.
\item{(b)}
$Lip^k(\R, \R)$ (the set of differentiable functions $f\colon \R \to \R$
of order $k$ whose derivatives are locally Lipschitz, 
denoted ${\cal M}_L^k$, for any $k\in \N \cup \{ \infty \}$.
\item{(c)}
$l^{\infty}(\N, \R)$, the set of bounded real sequences, denoted 
$l^{\infty}$. 

For any set $\cal M$, the class of $\cal M$-spaces as objets with the
$\cal M$-maps as morphisms forms a category denoted $\underline{{\cal M}}$.
This category has initial and final structures: if $(X_i)_{i\in I}$ is
a family of $\cal M$-spaces and $g_i\colon X \to X_i$ is given, the
initial structure on $X$ is generated by the set of functions
$\{ f\circ g_i\colon i\in I, f\in {\cal F}_i \}$.
The set of structure curves is
$\{ c\colon S \to X \colon g_i\circ c\in {\cal C}_i, \forall i\in I \}$.

Instead, if $g_i\colon X_i \to X$ is given, the final structure on $X$
is generated by the set of curves
$\{ g_i\circ c \colon i\in I, c\in {\cal C}_i \}$,
and the set of structure functions is
$\{ f\colon X \to R \colon f\circ g_i \in {\cal F}_i, \forall i\in I \}$.

Initial structures provides an $\cal M$-structure on products and subsets
and final structures on quotients.

Cartesian closeness provides a way to define $\cal M$-structures on spaces
of functions: let $\cal M$ containing all constant maps and let
$${\cal C}_{{\cal M}} = \{ c\colon S \to {\cal M} / 
\bar c\colon S\times S \to R \ \hbox{is an {\cal M}-map} \}$$
where $\bar c(s,t) = c(s)(t)$. If $\cal M$ satisfies that
$\Gamma \Phi {\cal C}_{\cal M} \subset {\cal C}_{\cal M}$, 
then for any $\cal M$-spaces $Y, Z$,
the set ${\cal M}(Y,Z)$ is an $\cal M$-space \cite\Frolicherb \
with structure curves
$${\cal C}_{{\cal M}(Y,Z)} = \{ h\colon S \to {\cal M}(Y,Z) /
\bar{h}\colon S\times Y \to Z \ \hbox{is an} \ {\cal M}-\hbox{map} \}$$

This condition is satisfied by ${\cal M}^{\infty}$ and $l^{\infty}$ 
and both are linearly generated \cite\Lawvere. An application of this
fact is 
that evaluation and composition maps are ${\cal M}^{\infty}$-maps
in this category.

In addition, if $E,F$ are normed spaces, the set of continuous linear
maps $E^{\prime}, F^{\prime}$ generates ${\cal M}$-structures on $E, F$
for any $\cal M$ and in particular if we take ${\cal M}= Lip^k(\R, \R)$ and
$U$ an open subset of $E$, then a map $g\colon U \to F$ is $Lip^k$ (in
the classical sense) if and only if $g\in {\cal M}(U, F)$: the composition
with the structure curves on $U$ and structure functions on $F$ are members
of $\cal M$ \cite\Fro.
This result shows that the theory of $\cal M$-spaces provides a 
generalization of differentiation theory in Banach spaces (in the
case of $Lip^k$ and smooth maps). 

Cartesian
closeness provides a simple application to the theory of 
groups of diffeomorphisms:
if $X$ is any ${\cal M}^{\infty}$-space,
then $C^{\infty}(X, X)$ is an ${\cal M}^{\infty}$-space 
by cartesian closeness and the two maps 
$$i,j\colon Diff(X) \to C^{\infty}(X, X)$$
given by $i(f) = f, j(f) = f^{-1}$ provides an ${\cal M}^{\infty}$-structure 
on $Diff(X)$ (the initial one), 
which makes it a smooth group (composition and invertion are
${\cal M}^{\infty}$-maps).

Even though we may define ${\cal M}^k$-spaces with 
$k< \infty$, they do not provide a generalization of 
$C^k$-differentiation theory, since there are maps
$g\colon E\to F$ between Banach spaces which are not $C^k$ but
$g\in {\cal M}^k(E, F)$, even for finite dimensional spaces.

\section 4. Diffeological spaces.

The category of diffeological spaces \cite\SouriauAlgoritme, 
\cite\DonatoGeometrie, \cite\DonatoVirasoro,
contains smooth manifolds, finite and infinite dimensional, and also 
finite dimensional $C^k$ manifolds as subcategories. 
They were defined by Souriau in order to extend quantization
through coadjoint orbits to groups of diffeomorphisms. In this
section we recall its definition.

The notation $P{\cal M}$ will mean a class of $C^k$ functions
$f\colon U\subset \R^n \to \R^m$ for some $0\leq k \leq \infty$,
and any $n,m\in \N$,
containing constants and closed under composition, such as:

\item{(a)}
The set of $C^k$ maps $f\colon U\subset \R^n \to \R^m$, for any
$U$ open and any $n, m\in \N$. This set is denoted $P{\cal M}^k$,
with $k\in \N \cup \{ \infty \}$.
\item{(b)}
The set of $C^k$ maps $f\colon U\subset \R^n \to \R^m$ whose
derivatives are locally Lipschitz, for any $U$ open and any 
$n, m\in \N$. Denoted $P{\cal M}^k_L$ with $k\in \N \cup \{ \infty \}$.

Given a set $X$, an $n$-{\it plaque} on $X$ is a function $p\colon U\to X$
where $U$ is an open subset of $\R^n$. A 
$P{\cal M}$-{\it diffeology} on $X$ is a
set $P(X)$ of $n$-plaques, for each $n\in \N$ such that the images of the
plaques covers $X$ and
\item{i.}
If a set $(p_i)$ of $n$-plaques admits a common extension, then the
smallest such extension is also an $n$-plaque in $P(X)$.
\item{ii.}
For each $\phi \in P{\cal M}$, $\phi \colon U^{\prime} \to U$, where
$U', U$ are open in $\R^m, \R^n$ respectively and 
for every  plaque $p\colon U \to X$, the map $p\circ \phi$ is also in $P(X)$.

The pair $(X, P(X))$ is called a $P{\cal M}$-{\it diffeological space}.
The set of $n$-plaques is denoted $P_n(X)$. 

The {\it standard} $P{\cal M}$-{\it diffeology} 
on $\R$ is defined as the set
of maps $f\colon U \subset \R^n \to \R$ in $P{\cal M}$ and it is
denoted by $P(\R)$.

Given two diffeological spaces $(X, P(X))$ and $(Y, P(Y))$,
a map $f\colon X \to Y$ is called {\it a $P{\cal M}$ morphism} 
if $p\in P(X)$ implies
$f\circ p\in P(Y)$. The set of such maps is denoted $D(P(X), P(Y))$ 
or $D(X,Y)$ if possible.
In particular if $Y= \R$ with the standard diffeology, 
then we use the notation
$D(X)$ instead. We shall say that two diffeologies $P^1(X), P^2(X)$
on $X$ are {\it equivalent} if $D(P^1(X)) = D(P^2(X))$,
and it is denoted by $(X, P_1(X)) \sim (X, P_2(X))$.

The class of $P{\cal M}$ diffeological spaces 
as objets with the differentiable
functions as morphisms forms a category, denoted $\underline{P{\cal M}}$.

\declare Remarks 4.1.

\item{[1]}
Any set $P_0$ of plaques on $X$, whose images covers $X$, generates a
diffeology $\bar{P_0}$, formed by plaques $p\colon U\to X$ such that
$\forall r\in U$ there exists an open neighborhood $U_r$, 
$\phi \in P{\cal M}$, and $p_0\in P_0$ such that
$p\vert_{U_r} = p_0\circ \phi$

\item{[2]}
Any set ${\cal F}_0$ of functions on $X$ generates a set $P$ of plaques 
$p$ on $X$ such that $f\circ p\in P{\cal M}$ for all $f\in {\cal F}$.
This set of plaques is denoted $\Gamma{\cal F}_0$. Also, any set
$P_0$ of plaques defines a set $\Phi{\cal F}$ of functions $f$ on
$X$ such that $f\circ p\in P{\cal M}$ for all $p\in P_0$.
Given $P_0$, then $(X, \bar{P_0})$ is a diffeology whose set of 
smooth functions is $\Phi \bar{P_0} = \Phi P_0$. Notice
that $\bar{P_0} = \Gamma \Phi P_0$. Also, notice that
given ${\cal F}_0$, we have that $(X, \Gamma {\cal F}_0)$ is a
diffeology whose set of differentiable functions is
$\Phi \Gamma{\cal F}_0$.

\item{[3]}
Given a collection $(X_j, P_j)_{j\in J}$ of diffeological spaces and maps
$g_j\colon X_j \to X$, the final diffeology on $X$ is generated by 
$\{ g_j\circ p / p\in P_j, j\in J \}$.
Instead, if $g_j \colon X \to X_j$, the initial diffeology on $X$ is the
set of plaques $p\colon U\subset \R^n \to X$ such that $g_j\circ p\in P_j$
for all $j\in J$. In particular this allows 
the definition of diffeologies on products,
quotients and subsets.

\item{[4]}
Cartesian closeness follows in general: if $(Y, P(Y)), (Z, P(Z))$ are
 diffeological spaces, the {\it functional diffeology} on $D(Y, Z)$
is formed by all plaques $p\colon U \to D(Y, Z)$ such that
$\bar p\colon U\times Y \to Z$ is differentiable, where
$\bar p(r,y) = p(r)(y)$. It follows that, given $(X, P(X))$, a
map $f\colon X \to D(Y, Z)$ is differentiable if and only if
$\bar f$ is differentiable, and therefore
$D(X, D(Y,Z))$ is diffeomorphic to $D(X\times Y, Z)$.

\item{[5]}
Given $(X, P(X)), (Y, P(Y))$ as before with corresponding
sets of smooth functions $D(X), D(Y)$, and a map
$f\colon X \to Y$, we might consider three natural definitions of a
differentiable function:
\itemitem{(i)}
$f$ is differentiable if $p\in P(X) \then f\circ p\in P(Y)$.
\itemitem{(ii)}
$f$ is differentiable if $g\in D(Y) \then g\circ f\in D(X)$.
\itemitem{(iii)}
$f$ is differentiable if $p\in P(X), g\in D(Y) \then g\circ f\circ p
\in P{\cal M}$.
\hfill\break
It is easy to check that $i\Rightarrow ii$ and $ii\Leftrightarrow iii$,
in addition, $ii\Rightarrow i$ if and only if $\Gamma \Phi P(Y) = P(Y)$.
This last condition may be added to the definition of a diffeology
$P(X)$ on $X$. In this case the
set $D(X)$ does not change, $P(X)$ becomes larger, and also the set of
differentiable functions $D(X, Y)$ for any $(Y, P(Y))$, but cartesian
closeness would not follow in general.

\item{[6]}
Notice that definition (i) corresponds to the category 
$\underline{P{\cal M}}$, definitions (ii) and (iii) will define
another category $\underline{P{\cal M}_1}$ with the same set $D(X)$
but with different set of morphisms $D_1(X,Y)$. The map 
$\gamma \colon \underline{P{\cal M}} \to \underline{P{\cal M}_1}$
given by the identity on morphisms and
$\gamma(X, P(X))= (X, \Gamma \Phi P(X))$ is a functor
such that 
$$\gamma(X, P(X))= \gamma(X, \Gamma \Phi P(X))$$ 
The
functor $\beta \colon \underline{P{\cal M}_1} \to \underline{P{\cal M}}$
given by the identity on objets and morphisms satisfies 
$\gamma \circ \beta = id$.
Final diffeologies are the same under both categories, and also the 
product diffeology. Cartesian closeness also follows in 
$\underline{P{\cal M}_1}$ if $P{\cal M}$ satisfies
$\Gamma \Phi (P{\cal M}) = P{\cal M}$.

\item{[7]}
If $(X, P^1(X)) \simeq (X, P^2(X))$ and $(Y, P(Y))$ satisfies
$\Gamma \Phi P(Y) = P(Y)$ ($P(Y)$ is maximal) then
$$D(P^1(X), Y) = D(P^2(X), Y)$$ 
If $P(Y)$ is not maximal then the
equality does not hold. For example, let $X=Y=\R^2$ and
$P^1(X) = P^2{\cal M}$, $P^2(X) = P(Y)= P^2_1{\cal M}\circ P^1{\cal M}$,
then $id\in D(P^1(X), Y)\setminus D(P^2(X), Y)$.

\section 5. Smooth structures and diffeologies.

\declare Lemma 5.1.(\cite\Fro).

Let $E, F$ be Banach spaces and
$k\in \N \cup \{ \infty \}$. A map $g\colon E \to F$ is
k-Lipschitz (in the classical sense) if and only if 
for every k-Lipschitz curve
$\gamma \colon \R \to E$ the map
$g\circ \gamma$ is a k-Lipschitz curve on $F$.

\declare Proposition 5.2.

Let $P{\cal M}$ be $P{\cal M}_L^k$
for some $k\in \N \cup \{ \infty \}$.
\item{a.}
Every $\cal M$-structure $({\cal C}, {\cal F})$ on $X$ defines a
$P{\cal M}$-diffeology $P(X)$ on $X$ with
$P_1(X) = {\cal C}$, $D(X) = {\cal F}$. This defines a functor 
$\Psi \colon \underline{{\cal M}} \to \underline{P{\cal M}}$ which
is 1-1 on objects and it is an embedding.
\item{b.}
Every $P{\cal M}$-diffeology $(X, P(X))$ defines an
$\cal M$-structure $({\cal C}, {\cal F})$ on $X$ generated by
 ${\cal C}_0 = P_1(X)$ such that
 ${\cal F}= D(X)$. This defines a faithful functor 
$\Upsilon \colon \underline{P{\cal M}} \to \underline{{\cal M}}$
    which is onto on 
objects and satisfies 
$\Upsilon \circ \Psi = Id$ and 
$\Psi \circ \Upsilon (X, P(X)) \sim (X, P(X))$.
\item{c.}
$\underline{{\cal M}} \simeq \underline{P{\cal M}_1}$ (and also for any
$P{\cal M}$ satisfying $*$, given below.

{\it Proof:} for any $P{\cal M}$ let us use the notation
$$P^n_m{\cal M} = \{ p\colon U\subset \R^m \to \R^n / p\in P{\cal M} \}$$
Observe that lemma 5.1. implies that $h\in P^1_m{\cal M}$ if and only if
$$h\circ \gamma \in {\cal M}, \forall \gamma \in P_1^m{\cal M} \eqno(*)$$
(a): Given $(X, {\cal C}, {\cal F})$ define
$$P(X) = \{ p\colon U\subset \R^n \to X / f\circ p\in P^1_m{\cal M},
\forall f\in {\cal F}, n\in \N \}$$
Then $P(X)$ is a diffeology satisfying $P_1(X) = {\cal C}$ and
${\cal F}\subset D(X)$. Conversely, if $f\in D(X)$ then 
$f\circ p\in P^1_m{\cal M}, \forall p\in P(X)$ and therefore
$f\circ p\in P^1_1{\cal M} = {\cal M}, \forall p\in P_1(X) = {\cal C}$.
this proves that $f\in {\cal F}$ and that ${\cal F} = D(X)$.


This defines
a map $\Psi \colon \underline{{\cal M}} \to \underline{P{\cal M}}$
which is 1-1 : if 
$$ \Psi (X, {\cal C}_1, {\cal F}_1) = \Psi (X, {\cal C}_2, {\cal F}_2)
= (X, P(X))$$
then ${\cal C}_1= P_1(X) = {\cal C}_2 $ and 
therefore ${\cal F}_1 = {\cal F}_2$.
Given $f\in {\cal M}(X, Y)$ define $\Psi (f) = f$. This is a map
which is 1-1 and
onto since $D(X, Y) = {\cal M}(X, Y)$:
 if $f\in {\cal M}(X, Y)$ and
$p\in P(X)$, then $h\circ f\in {\cal F}_X, \forall h\in {\cal F}_Y$). It
follows that
$h\circ f\circ p\in P^1_m{\cal M}$, $\forall h\in {\cal F}_Y$
this means that $f\circ p\in P(Y)$ and therefore that $f\in D(X, Y)$.
Conversly, if $f\in D(X, Y)$, then $f\circ p\in P(Y)$, $\forall p\in P(X)$;
thus 
$$f\circ c\in P_1(Y) = {\cal C}_Y, \forall c\in P_1(X) = {\cal C}_X$$
and this proves that $f\in {\cal M}(X, Y)$.

(b) Let $(X, P(X))$ be a $P{\cal M}$-diffeology and let 
$\Upsilon (X, P(X))$ be
the smooth $\cal M$-structure generated by ${\cal C}_0 = P_1(X)$.
Then 
$${\cal F} = \{ f\colon X \to \R / f\circ p\in {\cal M}, 
\forall p\in P_1(X) \}$$
It follows that $D(X) \subset {\cal F}$. Now let $f\in {\cal F}$. Let
$p\in P_m(X)$, then $p\circ \gamma \in P_1(X), \forall 
\gamma \in P_1^m{\cal M}$
, then
$(f\circ p)\circ \gamma \in {\cal M}, \forall \gamma \in P_1^m{\cal M}$,
and taking $h= f\circ p$ in $*$ we obtain 
$f\circ p\in P^1_m{\cal M}$ and 
then $f\in D(X)$.
This proves that $D(X)= {\cal F}$. If $(Y, P(Y))$ is another
diffeology, we define $\Upsilon (f) = f$ if $f\in D(X, Y)$. We have
that $D(X, Y)\subset {\cal M}(X, Y)$)and therefore $\Upsilon$ is also a
faithful functor: if $f\in D(X, Y)$ and $h\in {\cal F}_Y$
then $f\circ p\in P_1(Y), \forall p\in P_1(X)$, 
then $(h\circ f)\circ p\in {\cal M}, \forall p\in P_1(X)$,
consecuently $(h\circ f)\in {\cal F}_X$.

Given $(X, {\cal C}, {\cal F})$, let $\Psi (X, {\cal C}, {\cal F})=
(X, P(X))$ and 
$$\Upsilon (X, P(X))=(X, \bar{{\cal C}}, \bar{{\cal F}})$$
Then ${\cal F} = D(X) = \bar{{\cal F}}$ and therefore
${\cal C} = \bar{{\cal C}}$. This proves that $\Upsilon$ is onto and
$\Upsilon \circ \Psi = Id$. Conversely, given $(X, P(X))$, let
$\Upsilon (X, P(X)) = (X, {\cal C}, {\cal F})$
and $\Psi (X, {\cal C}, {\cal F}) = (X, \bar{P}(X))$, then 
$D_P(X) = {\cal F} = D_{\bar P}(X)$
and therefore $\Psi \circ \Upsilon (X, P(X)) \sim (X, P(X))$.

(c) The proof is also straightforward.
\qed

\declare Remarks 5.3.

\item{[1]}
The functor $\Psi$ is not onto on objects; for example, take
$X= \R^2$ with
$P^1_1(\R^2) = P_1^2{\cal M}$
and $P^1_n(\R^2) = P^1_1(\R^2) \circ P_n^1{\cal M}$. Then the image of
every plaque is a curve. If 
$\Psi (\R^2, {\cal C}, {\cal F}) = (\R^2,  P^2(\R^2))$ then
${\cal C} = P^1_1(\R^2)$ and therefore $P^2_n(\R^2) = P^2_n{\cal M}$,
which is different from $P^1_n(\R^2)$. However, given $(X, P(X))$ we
have 
$$\Phi (X, P_1(X), D(X)) \sim (X, P(X))$$
Notice that if the space $P{\cal M}$ satisfies $*$ then a smooth structure
is identical with an equivalence class of diffeological spaces. Otherwise
 there are two types of diffeological spaces: those generated, or not
by curves $P_1(X)$. In the first place 
they are equivalent to the smooth structure
$(X, P_1(X), D(X))$.

\item{[2]}
The functor $\Upsilon$ is not 1-1 on objects: both  
$P^i(\R^2), i=1,2$ given above have the same image 
$(\R^2, {\cal C}, {\cal F})$.

\item{[3]}
If we consider $P{\cal M}^k$ with $k\neq \infty$, then (a) remains
valid in the proposition. In part (b), $\Upsilon$ is still a faithful
functor, ${\cal C}$ is generated by $P_1(X)$ and $D(X) \subset {\cal F}$
(only), $\Upsilon$ is onto on objects, but $\Upsilon \circ \Psi \neq Id$,
$\Psi \circ \Upsilon \neq Id$. Also $\underline{{\cal M}}$ is
embedded in $\underline{P{\cal M}_1}$ but the last one is larger.

\item{[4]}
${\cal M}^k_L$-structures provides a generalization
of ${\cal M}^k_L$-Banach manifolds, but do not generalize ${\cal M}^k$-
manifolds, not even in the finite dimensional case. Diffeological spaces
generalizes ${\cal M}^k_L$-Banach manifolds and finite dimensional
${\cal M}^k$-manifolds.

\item{[5]}
Some disadvantages of those diffeologies outside 
$\underline{P{\cal M}_1}$ are: $P(X)$ determines $D(X)$ but $D(X)$ does
not determine $P(X)$, but instead, it determines a class of diffeologies.
Also, consider the following example: $P{\cal M} = P{\cal M}_L^k$ and
$X=Y=\R^2$, $P^1(X)= P^1(Y)=P{\cal M}^k$ and $P^2(X)=P^2(Y)=P{\cal M}^k_L$.
Then $D(X)= D(Y)=P_1^2{\cal M}^k_L$ (because a function
$f\colon \R^2 \to \R$ is $Lip^k$ iff it is $Lip^k$ along smooth curves
iff it is $Lip^k$ along $Lip^k$ curves. Therefore $D(P^2(X), P^2(Y))$
is $P_2^2{\cal M}^k_L$ but $D(P^1(X), P^1(Y)= P_2^2{\cal M}^{\infty}$
(by Boman's result \cite\Boman ). This is unexpected because we probably
hope to obtain spaces of $Lip^k$ functions.

\item{[6]}
There are important examples (such as $P{\cal M}^k$) for which $*$ does
not hold and $\underline{P{\cal M}_1}$ is larger than $\underline{\cal M}$.

\section 6. Tangent structures.

Both, smooth structures and diffeologies provides a relation between
the spaces of curves (plaques) and spaces of functions. 
In the 
first case, convenient spaces are the most general vector spaces, with
a smooth structure which is linearly generated for which the structure
curves have derivatives (the definitions are recalled in the next section).
In the second case, tangent diffeologies (with a tangent structure)
provides n-dimensional directions
at every point, allowing the definition of tangent spaces, derivatives of
functions, vector fields and differential forms, in which the tangent bundle
itself is a tangent diffeology.
The category of tangent diffeologies \cite\Torre, \cite\TorreTangen,
\cite\TorreBan provides a generalization of smooth and $C^k$ manifolds. 
The definition is recalled in this section. 
For simplicity of notation we may assume 
that all plaques $p$ through a point $F\in X$ are defined in a neighborhood
of $0$ and $p(0)=F$ and $p\circ \phi$ means that 
$\phi \colon U'\subset \R^m \to U \subset \R^n$ and $p\colon U\to X$. 

\item{(a)}
Given a $P{\cal M}$-diffeology $(X, P(X))$, a {\it tangent structure}
on $X$ is a collection $\sim$ of equivalence relations $\sim_F^n$
on $P^n_F(X)$, for each $F\in X$ and $n\leq k$ satisfying the
two consistency conditions:
\itemitem{(i)}
$p_1 \circ \phi \sim_F^m p_2 \circ \phi$ whenever $p_1 \sim_F^n p_2$
and $\phi \in P_n^m{\cal M}$.
\itemitem{(ii)}
$p_1 \sim_F p_1|_V$ if
$p_1 \: U \to X$ and $V \subset U$ is an open neighborhood of~$0$.
\item{(b)}
The tangent structure is called {\it linear\/} if the set
$P_F^n(X) /{\sim_F^n}$ carries a vector space structure for
each $F\in X$ with the following consistency conditions: 
\itemitem{(i)}
If
$p_{12} \in [p_1] + [p_2]$ and $\phi \in P{\cal M}$ 
then
$p_{12} \circ \phi \in [p_1 \circ \phi] + [p_2 \circ \phi]$.
\itemitem{(ii)}
If $p\in c[p_1]$ then $p\circ \phi \in c[p_1\circ \phi ]$, 
for all $\phi \in P{\cal M}$.
\item{(c)}
The linear structure is called {\it continuous\/} if given two
$(n+m)$-plaques $p_i(r,s)$ with $i = 1,2$ such that
$p_1(r,0) = p_2(r,0)$ for each $r$,
then there exists a plaque $p_{12}(r,s)$ such that
$$
[p_{12}(r,s)]_s = [p_1(r,s)]_s + [p_2(r,s)]_s  \qquad\hbox{for all $r$}.
$$

The standard $P{\cal M}$ 
tangent structure on $(\R, P(\R))$ is defined by
$p_1\sim^n_F p_2$ if $D^ip_1(0)= D^ip_2(0)$ for all $i\leq n$.

We shall say also that $(X, P(X), \sim)$ is a 
$P{\cal M}$ (linear) tangent diffeological space (\SDS).
The class of the plaque $p(t)$ at a point $F \in X$ will be denoted by
$[p]$ or by $[p]_t$.
The set $\Tn_F^n X := P_F^n(X) /{\sim_F^n}$ is called the $n$-th
{\it tangent space\/} at~$F$, and the disjoint union
$\Tn^n X := \bigsqcup_{F\in X} \Tn_F^n X$ is called the $n$-th
{\it tangent bundle\/} over~$X$.

Let $(X,P(X),\sim)$, $(Y,P(Y),\approx)$ be two tangent diffeologies; a
 $P{\cal M}$ morphism ~$f \: X \to Y$ 
is called $P{\cal M}$-differentiable at~$F$
if for all $n \leq k$, $f \circ p_1 \approx_{f(F)}^n f \circ p_2$
whenever $p_1 \sim_F^n p_2$. The set of differentiable functions is
denoted ${\cal C}^k(X,Y)$ or ${\cal C}(X)$ if
$Y= \R$, with the standard tangent diffeology (it is denoted by
$C^k(P(X), P(Y))$ if necessary).
If the tangent diffeologies are linear, $f$ is called a $C^k$-map 
(or a{\it smooth
map\/} if $k=\infty$) if, in addition,
$$
D_F^n f : \Tn_F^n X \to \Tn_{f(F)}^n Y : [p] \mapsto [f \circ p]
$$
is linear for each $F \in X$. The set of $C^k$-maps is denoted by
${\cal C}^k_l(X,Y)$.

\declare Remark 6.1.

\item{[1]}
Notice that ${\cal C}^k(X, Y) \subset D(X,Y)$.

\item{[2]}
The set of smooth diffeological spaces as objets with the
differentiable maps as morphisms is a category, since,
if $f\colon A \to B$ and $g\colon B\to C$ are differentiable then
$g\circ f$ is a $P{\cal M}$ morphism and preserves the relations. In
addition
$$
D^n_F(g\circ f)[p] = (D^n_{f(F)}g)(D^n_Ff[p])
$$
(and this is linear), therefore $g\circ f$ is $C^k$ (a $C^k$-map).

\item{[3]}
Let $(X, P(X))$ be a diffeological space.
If $p$ is an n-plaque at $F$ and $f\in D(X)$, we shall denote by
$D^n_pf$ the sequence $D^i(f\circ p)\vert_0, i= 1, \cdots n$.
Let $\cal F$ a subspace of $D(X)$. The {\it canonical \SDS}\ 
defined by $\cal F$ is defined by:
$$p_1\sim_F^n p_2 \iff D^n_{p_1}f = D^n_{p_2}f, \forall f\in 
{\cal F}$$
and $[p_1]+[p_2]$ is defined by:
$p_{12}\in [p_1]+[p_2]$  if and only if 
$D_{p_{12}}f = D_{p_1}f+ D_{p_2}f$, for all $f\in {\cal F}$
and  if $c\in \R$, define $c[p_1(t)] = [p_1(ct)]$. The set of these 
classes is denoted $\sim_{{\cal F}}$. In particular if
${\cal F} = D(X)$ then $D(X)= C^k(X)$.
(called {\it the canonical} \SDS).
(This is an equivalence relation that satisfies the two consistency 
conditions. Then $(X, P(X), \sim_{\cal F})$ is an 
\SDS \ with a partial vector space
structure).
Since $c[p]$ always exists, any diffeology has a
canonical \SDS \ for which each tangent space is a cone, and
satisfies $D_{p(ct)}f = cD_pf$ for every $c, p, f$.
However, not every canonical \SDS \ is linear, for example,
the intersection of two transverse smooth curves in $\R^n$ 
is not linear: at the point of
intersection the tangent space is the sum of two lines, everywhere else
it is a line. 

\item{[4]}
Notice also that if $P(X)$ satisfies  the condition:
$\forall p_1, p_2$ there exists $ p_{12}$ such that 
$D_{p_{12}}f= D_{p_1}f + D_{p_2}f, \forall f\in {\cal F}$ 
 then ${\cal T}_F^nX$ is
a vector space and the structure is linear. In particular this occurrs if
$P(X)$ generates all derivations
$\phi_F\colon D(X) \to \R$, that is, 
there exists $ p\in P(X)$ such that 
$$\phi_F(f) = D(f\circ p)\vert_0, \forall f \in D(X)$$

\section 7. Convenient spaces.

In this section we recall  
the definition of a convenient space and some
of the most important properties of these spaces. 
They are the most 
general vector spaces with a smooth structure 
in which $Lip^k$-curves have unique derivatives,
forming
subcategory of $\underline{{\cal M}_L^k}$. (They are the most convenient
spaces for differentiation in many ways).

A dual vector space is a pair $(X, X^{\prime})$ where $X$ is a vector
space and $X^{\prime}$ is a subspace of the algebraic dual. Each dual
vector space determines several structures: (a) the {\it weak topology}
$\tau_W$
is the initial one induced by $X^{\prime}$ (it is the coarsest locally
convex topology yielding $X^{\prime}$ as topological dual),
(b) the {\it Mackey topology} $\tau_M$ is the finest locally 
convex topology yielding
$X^{\prime}$ as topological dual \cite\Jarchow, 
(c) $X^{\prime}$ generates a $Lip^k$-structure,
(d) the {\it Mackey closure topology} $\tau_{MC}$ is the final topology 
generated by the set of $Lip^k$-curves. In addition it determines
a bornology, a $l^{\infty}$ vector space and a convergence 
structure. A good description of the functors that relate these 
structures is given in \cite\Fro. 

A {\it preconvenient space} is a dual vector space $(X, X^{\prime})$ 
such that
the $Lip^k$ structure that $X^{\prime}$ determines (linearly generated)
has exactly $X^{\prime}$ as set of linear structure functions.
An equivalent condition is that $X^{\prime}$ is the topological
dual under $\tau_{MC}$. 
The category of preconvenient spaces is denoted
$\underline{Pre}$.

\declare Remarks 7.1.

\item{[1]}
Every preconvenient space have the three locally convex
topologies defined above and each of
these topologies have $X^{\prime}$ as the space of continuous linear
functionals. Every l.c.s. $X$ determines a dual pair 
with $X^{\prime}$ the topological dual, 
however, not every l.c.s. forms
a preconvenient space with the topological dual. For example,
an uncountable direct sum of $\R$ with the box topology is not pre-convenient
\cite\Jarchow. However, every metrizable space is pre-convenient and
in this case the Mackey topology and the Mackey closure topology coincide.
(in general $\tau_W \subset \tau_M \subset \tau_{MC}$). A dual vector space
is preconvenient iff the Mackey topology is an inductive limit of
seminormable spaces.

\item{[2]}
If a d.v.s. $(X, X^{\prime})$ is preconvenient, the space determines and
is determined by each of the structures mentioned above:
a topological convex vector bornology 
($B\subset X$ is bounded
iff $l(X)$ is bounded $\forall l\in X^{\prime}$), a bornological 
locally convex
topology $\tau_M$,
a linearly generated $l^{\infty}$ structure (generated by $X^{\prime}$),
a linearly generated $Lip^k$ structure (generated by $X^{\prime}$),
a convergence structure (a sequence $(x_n)$ converges to $0$ if there
exists reals $t_n\to \infty$ such that $\{ t_na_n/n\in \N \}$ is bounded),
an arc-determined topology (the Mackey closure topology). The space 
$X^{\prime}$ is recovered by taking the set of linear functions respecting
the structures.
Preconvenient spaces form the category 
that is identified with subcategories of these 
categories, for example, $\underline{Pre}$ is isomorphic to the category of
linearly generated $Lip^k$ vector spaces.
\item{[3]}
Given $X, Y$ pre-convenient, a map $f\colon X \to Y$ is a $Lip^k$ map ( or a 
$\underline{Pre}$ morphism) iff
it is continuous for the Mackey closure topology \cite\Kri, iff
$f^*(Y^{\prime})\subset X^{\prime}$, iff it is continuous for the Mackey
topology.

\item{[4]}
Given a vector space $X$ and $X_j\in \underline{Pre}, \forall j\in I$
and linear maps $m_j\colon X \to X_j$ (respectively $m_j\colon X_j\to X$),
the initial (resp. final) $Lip^k$-structure on $X$ defines an initial
(resp. final) preconvenient structure on $X$.

Given  $(X, X^{\prime})\in \underline{Pre}$ and $c\colon \R \to X$. We
say that $c^{\prime}(t)\in X$ is the {\it weak derivative}
 of $c$ at $t$ if $(l\circ c)^{\prime}(t)$ exists and equals 
$l(c^{\prime}(t)), \forall l\in X^{\prime}$: $c^{\prime}(t)$ is the
derivative of $c$ with respect to $\tau_W$. 
Also, $\int_t^sc$ is the {\it weak integral} of $c$
if $\int_t^s(l\circ c)(u)du$ exists and equals 
$l(\int_t^sc), \forall l\in X^{\prime}$.

A preconvenient space is called {\it separated} if every $Lip^k$-curve
has at most one weak derivative. The subcategory of these spaces is
denoted $\underline{Spre}$.

A space $(X, X^{\prime})\in \underline{Spre}$ is called
{\it complete} if every $Lip^k$-curve has a weak derivative.
The subcategory of them is denoted $\underline{Con}$.

\declare Remarks 7.2.

\item{[1]}
Given $(X, X^{\prime}) \in \underline{Pre}$, the following are equivalent:
\itemitem{a.}
$X^{\prime}$ separates points.
\itemitem{b.}
Every $Lip^k$-curve has at most one weak derivative.
\itemitem{c.}
Every $Lip^k$-curve has at most one weak integral.
\itemitem{d.}
$\tau_{MC}$ is Hausdorff.
\itemitem{e.}
$\tau_M$ is Hausdorff.

\item{[2]}
Given $(X, X^{\prime})\in \underline{Spre}$, the following are
equivalent:
\itemitem{a.}
$(X, X^{\prime})$ is convenient.
\itemitem{b.}
For every $Lip^k$-curve $c$, the weak integral $\int_0^sc$ exists.
\itemitem{c.}
If $\{ l(x_n) / n\in \N \}$ is bounded $\forall l\in X^{\prime}$ and
$(t_n)\in l^1$, then $\sum t_nx_n$ converges weakly.
\itemitem{d.}
Every Mackey Cauchy sequence converges (a sequence $(x_n)$ {\it converges} to
$x_0$ if there exists $t_n\in \R^*, \forall n\in \N$ such that
$\lim_{n\to \infty}t_n = \infty$ and
  $\{ l(t_n(x_n-x_0)) / n\in \N \}$ is bounded $\forall l\in X^{\prime}$.
The sequence is called a {\it Mackey Cauchy sequence} if there exists
$t_{n,m}\in \R^*, \forall n,m\in \N$ with
$\lim_{n,m\to {\infty}}t_{n,m}= \infty$ and the set
$\{ l(t_{n,m}(x_n-x_m)) / n,m\in \N \}$
is bounded $\forall l\in X^{\prime}$).

\item{[3]}
The product $\pi X_j$ of convenient spaces $X_j$ is convenient (the 
$Lip^k$ curves on $X$ are those whose coordinates are $Lip^k$).
The coproduct $\coprod X_j$ is also convenient (the vector space is the direct 
sum, the dual is the product of the duals, the structure curves 
$c\colon \R \to X$ are those which are locally representable as finite
sums of $Lip^k$curves $c_j\colon \R \to X_j$). 
Inductive limits may also be convenient: if $X_j$ is convenient, and it is a
pre-subspace of $X_{j+1}$, closed with $\tau_{MC}$ , then the
inductive limit in $\underline{Pre}$ is convenient (the $Lip^k$
curves are locally $Lip^k$ curves into some $X_j$). The tensor product
of convenient spaces is also convenient.

\item{[4]}
Subspaces are convenient in some cases: If $V$ is a vector space,
$X$ is convenient and $f\colon V \to X$ is injective, linear and $f(V)$
is closed with the Mackey closure topology, then the initial 
pre-structure of $f(V)$ is convenient. Quotients may also be formed in
some cases: if $f\colon X \to V$ where $X$ is convenient and $V$ is
a vector space, then the final pre-structure on $V$ is convenient iff
the kernel of $f$ is closed in the Mackey topology and the final locally
convex topology of $V$ is Mackey complete.

\item{[5]}
If $X$ is convenient and $Y$ is preconvenient then $L(X,Y)$ is a convenient
space with the initial structure induced by the evaluations. More generally
$L^n(X,Y)$ is convenient and 
$$L(X, L^{n-1}(X, Y)) = L^n(X, Y)$$

\item{[6]}
If $X$ is a convenient space and $k\in \N \cup \{ \infty \}$, then
$Lip^k(\R, X)$ is also a convenient space: the $\underline{Pre}$-
structure is the initial one induced by $\delta^0, \cdots , \delta^{k+1}$
where $\delta^k\colon Lip^k(\R, X) \to l^{\infty}(\R^{<p>},X)$,
(and $l^{\infty}(\R^{<p>},X)$ has the initial structure defined by
evaluations)
is defined as follows: if $f\colon D\subset \R \to X$, let
$$D^{<k>}= \{ (t_0, \cdots ,t_k)\in D^{k+1}: 
t_j\neq t_i \ \hbox{if} \ i\neq j \}$$
define $\delta^0f=f$ and
$$\delta^kf(t_0, \cdots ,t_k)= 
{k\over t_0-t_k}(\delta^{k-1}f(t_0, \cdots ,t_{k-1})
-\delta^{k-1}f(t_1, \cdots ,t_k))$$
Given $k\in \N \cup \{ \infty \}$, $(Y, {\cal C}, {\cal F})$ a
$Lip^k$-space and $X$ a convenient space, then $Lip^k(Y, X)$ is also
a convenient space, where $Lip^k(Y, X)$ denotes the vector space formed
by the $Lip^k$-maps $Y\to X$ together with the initial $\underline{Pre}$
-structure induced by the linear maps 
$c^*\colon Lip^k(Y, X) \to Lip^k(\R, Y), \ \hbox{for}\  c\in {\cal C}$.

Differentials of functions in preconvenient spaces are defined as
follows:
let $U\subset \R$ be open, $(X, X^{\prime}) \in \underline{Pre}$. 
A curve $c\colon U\to X$ is called {\it differentiable} if 
the weak derivative of $c$ at $t$ exists $\forall t\in U$. And it is
called {\it k-times differentiable} if $c$ is differentiable and
$c^{\prime}(t)$ is $(k-1)$-times differentiable. In particular if $X$
is a convenient space 
then $c$ is a $Lip^k$-curve iff $c$ is $j$-times
differentiable and $c^j$ is a $Lip^{k-j}$ curve if $0\leq j\leq k$.

The differential of a function $f\colon U\subset X \to Y$, where 
$(Y, Y^{\prime})$ is also preconvenient is defined as follows: $f$ is
called differentiable at $x$ in direction $v\in X$ if $c(t)= f(x+tv)$
is differentiable. This derivative is denoted $df(x,v)$ and $f$ is
called {\it differentiable} if $df(x,v)$ exists $\forall x\in U$, 
$\forall v\in X$, (if
$df$ is a $Lip^0$-map also, then $f$ is called {\it Lip-differentiable}).
The {\it differential} of $f$ is the map
$df\colon U\times X \to Y$. And $f$ is called {\it $k$-times differentiable}
if $f$ is differentiable and $df(\cdot ,v)$ is $(k-1)$-times 
differentiable $\forall v\in X$, (if $f$ is Lip-differentiable and
$df$ is $(k-1)$-times Lip-differentiable $\forall v\in X$ then $f$ is
called {\it $k$-times Lip-differentiable}). In particular if $X, Y$ are
convenient spaces then these two definitions and the definition of
a $Lip^k$ map are equivalent. It is also equivalent to a stronger
condition of differentiability in which the limit exists uniformly
with respect to $x$ and $v$ and $f^{\prime}(x)\in L(X, Y)$ where
$f^{\prime}(x)(v)= df(x,v)$ \cite\Seip (with respect to b-compact
sets: $K$ is b-compact in $X$ if there exists $B$ absolutly convex
such that $K$ is compact in
$X_B = \cup_{n\in \N}nB$).

\bye